\documentclass{article}%
\usepackage{amsmath}
\usepackage{amsfonts}%
\usepackage{amssymb}
\newtheorem{theorem}{Theorem}

\newtheorem{corollary}[theorem]{Corollary}

\newtheorem{definition}[theorem]{Definition}
\newtheorem{example}[theorem]{Example}

\newtheorem{lemma}[theorem]{Lemma}

\newtheorem{proposition}[theorem]{Proposition}
\newtheorem{remark}[theorem]{Remark}

\newenvironment{proof}[1][Proof]{\textbf{#1.} }{\ \rule{0.5em}{0.5em}}

\begin{document}

\title{Geodesic Webs on a Two-Dimensional Manifold and Euler Equations}
\author{Vladislav V. Goldberg\\New Jersey Institute of Technology, Newark, NJ, USA
\and Valentin V. Lychagin\\University of Tromsoe, Tromsoe, Norway}
\date{}
\maketitle

\begin{abstract}
We prove that any planar $4$-web defines a unique projective
structure in the plane in such a way that the leaves of the
foliations are geodesics of this projective structure. We also
find conditions for the projective structure mentioned above to
contain an affine symmetric connection, and conditions for a
planar $4$-web to be equivalent to  a geodesic $4$-web on an
affine symmetric surface. Similar results are obtained for planar
$d$-webs, $d > 4$, provided that additional $d-4$ second-order
invariants vanish.

\end{abstract}

\section{Introduction}
In this paper, which is a continuation of the paper \cite{AGL04},
 we study geodesic webs, i.e., webs whose leaves are
totally geodesic in a torsion-free affine connection.

We study in detail the planar case and prove that any planar
$4$-web defines a unique projective structure in the plane in such
a way that the leaves of the foliations are geodesics of this
projective structure. We also find conditions for the projective
structure mentioned above to contain an affine symmetric
connection, and conditions for a planar $4$-web to be equivalent
to  a geodesic $4$-web on an affine symmetric surface.

Similar results are obtained for planar $d$-webs, $d > 4$,
provided that additional $d-4$ second-order invariants vanish.

We also apply the obtained results to a surface of constant
curvature and to the linear webs. This allows us to prove the
Gronwall-type theorem (see \cite{Gr12} and \cite{GL06}) and its
natural generalizations in the case of geodesic webs.

For this, first, we find necessary and sufficient conditions for
the foliation defined by level sets of a function
 to be totally geodesic in the torsion-free connection. This
brings us to what we call the flex equation. The flex equation
possesses the infinite-dimensional pseudogroup of gauge
symmetries. Factorization of the flex equation with respect to
this pseudogroup leads us to the Euler equation as well as to
natural generalizations of it. This reduction gives us a way to
solve the flex equation.

Second, we apply these conditions to find a linearity criterion
for planar webs mentioned above. We formulate these conditions by
means of the flex equation(s) and show how to describe linear webs
in terms of the Euler equation.

For all these webs we find conditions which web function of a
linear web must satisfy.

\section{Linear Connections in Nonholonomic Coordinates}

Let $M$ be a smooth manifold of dimension $n.$ Let vector fields
$\partial _{1},...,\partial_{n}$ form a basis in the tangent
bundle and let $\omega ^{1},..,\omega^{n}$ be the dual basis. Then
\[
\lbrack\partial_{i},\partial_{j}]=\sum_{k}c_{ij}^{k}\partial_{k}%
\]
for some functions $c_{ij}^{k}\in C^{\infty}\left(  M\right)  ,$
and
\[
d\omega^{k}+\sum_{i<j}c_{ij}^{k}\omega^{i}\wedge\omega^{j}=0.
\]
Let $\nabla$ be a linear connection in the tangent bundle, and let
$\Gamma_{ij}^{k}$ be the Christoffel symbols of second type. Then%
\[
\nabla_{i}\left(  \partial_{j}\right)
=\sum_{k}\Gamma_{ji}^{k}\partial_{k},
\]
where $\nabla_{i}\overset{\text{def}}{=}\nabla_{\partial_{i}},$
and
\[
\nabla_{i}\left(  \omega^{k}\right)
=-\sum_{j}\Gamma_{ij}^{k}\omega^{j}.
\]
The covariant differential of a vector field
\[
d_{\nabla}:\mathbf{D}\left(  M\right)  \rightarrow\mathbf{D}\left(
M\right) \otimes\Omega^{1}\left(  M\right)  ,
\]
and the covariant differential of a differential 1-form
\[
d_{\nabla}:\Omega^{1}\left(  M\right)  \rightarrow\Omega^{1}\left(
M\right) \otimes\Omega^{1}\left(  M\right)
\]
take the following form:
\[
d_{\nabla}\left(  \partial_{i}\right)
=\sum_{k,j}\Gamma_{ij}^{k}~\partial _{k}\otimes\omega^{j},\] and
\[
d_{\nabla}\left(  \omega^{k}\right) =-\sum
_{i,j}\Gamma_{ij}^{k}\omega^{j}\otimes\omega^{i}.
\]
Remark that if the connection $\nabla$ is torsion-free, then
\[
\Gamma_{ij}^{k}-\Gamma_{ji}^{k}=c_{ij}^{k}.
\]

For the curvature tensor $R=d_{\nabla}^{2}:\mathbf{D}\left(
M\right) \rightarrow\mathbf{D}\left(  M\right)
\otimes\Omega^{2}\left(  M\right)  $
one has%
\[
R\left(  \partial_{i},\partial_{j}\right)
:\partial_{k}\longmapsto\sum _{l}R_{kij}^{l}\partial_{l},
\]
where
\[
R_{kij}^{l}=\partial_{i}\left(  \Gamma_{jk}^{l}\right)
-\partial_{j}\left(
\Gamma_{ik}^{l}\right)  +\sum_{m}\left(  \Gamma_{jk}^{m}\Gamma_{im}^{l}%
-\Gamma_{ik}^{m}\Gamma_{jm}^{l}-c_{ij}^{m}\Gamma_{mk}^{l}\right) .
\]

\section{Geodesic Foliations and Flex Equations}

The covariant differential%
\[
d_{\nabla}:\Omega^{1}\left(  M\right)  \rightarrow\Omega^{1}\left(
M\right) \otimes\Omega^{1}\left(  M\right)
\]
splits into the direct sum $d_{\nabla}=d_{\nabla}^{a}\oplus
d_{\nabla}^{s}$ according to the splitting of tensors into the sum
of skew-symmetric and
symmetric ones:%
\[
\Omega^{1}\left(  M\right)  \otimes\Omega^{1}\left(  M\right)
=\Omega ^{2}\left(  M\right)  \oplus S^{2}\left(  M\right)  .
\]
Then the connection is a torsion-free if and only if the
skew-symmetric component coincides (up to sign) with the de Rham
differential:
\[
d_{\nabla}^{a}=-d:\Omega^{1}\left(  M\right)
\rightarrow\Omega^{2}\left( M\right)  .
\]
In the nonholonomic coordinates the symmetric component
$d_{\nabla}^{s}$ has the form
\[
d_{\nabla}^{s}\left(  \omega^{k}\right)
=-\sum_{i,j}\Gamma_{ij}^{k}\omega
^{i}\cdot\omega^{j}%
\]
where $\cdot$ means the symmetric product of differential
$1$-forms.

\begin{lemma}
The foliation defined by a differential $1$-form $\theta$ is
totally geodesic in the connection $\nabla$ if and only if
\begin{align*}
d_{\nabla}^{a}\left(  \theta\right)  =\alpha\wedge\theta,\;\; d_{\nabla}%
^{s}\left(  \theta\right)  =\beta\cdot\theta
\end{align*}
for some differential $1$-forms $\alpha$ and $\beta.$
\end{lemma}

\begin{proof}
See Proposition 2 in \cite{AGL04}.
\end{proof}

\begin{corollary}
The foliation defined by a differential $1$-form $\theta$ is
totally geodesic
in the torsion-free connection $\nabla$ if and only if%
\[
d_{\nabla}^{s}\left(  \theta\right)  =\beta\cdot\theta
\]
for a differential $1$-form $\beta.$
\end{corollary}

\begin{corollary}
The foliation defined by level sets of a function $f$ is totally
geodesic in the torsion-free connection $\nabla$ if and only if
the quadratic form
\[
d_{\nabla}^{s}d\left(  f\right)  \in S^{2}\left(  M\right)
\]
vanishes on the level sets $f=\operatorname*{const}.$
\end{corollary}

\begin{lemma}
Assume that differential $1$-forms $\tau_{1},..,\tau_{n}$ are
linearly independent on $M$. A quadratic form
\[
Q=\sum_{ij}Q_{ij}\tau_{i}\cdot\tau_{j}%
\]
vanishes on the distribution $\tau_{1}+\cdots+\tau_{n}=0 $ if and
only if we have
\[
Q_{ii}+Q_{jj}=2Q_{ij}%
\]
for all $i,j.$
\end{lemma}

\begin{proof}
One has
\[
\sum_{ij}Q_{ij}\tau_{i}\cdot\tau_{j}=\left(
\tau_{1}+\cdots+\tau_{n}\right) \cdot\left(
x_{1}\tau_{1}+\cdots+x_{n}\tau_{n}\right)
\]
for some functions $x_{1},...,x_{n}.$

Then $2Q_{ij}=x_{i}+x_{j}.$ Taking $i=j$, we get $x_{i}=Q_{ii}$,
and then $Q_{ii}+Q_{jj}=2Q_{ij}.$
\end{proof}

Now we apply this lemma to the quadratic form
$Q=d_{\nabla}^{s}d\left( f\right)  $ and
$\tau_{i}=f_{i}\omega^{i},$ where $f_{i}=\partial_{i}\left(
f\right)  .$ We have
\[
df=\sum_{i}f_{i}\omega^{i}=\sum_{i}\tau_{i}%
\]
and
\begin{align*}
d_{\nabla}^{s}\left(  df\right)   &  =d_{\nabla}^{s}\left(  \sum_{k}%
f_{k}\omega^{k}\right)  =\sum_{k}f_{k}d_{\nabla}^{s}\left(
\omega^{k}\right) +\sum_{i}df_{i}\cdot\omega^{i}\\
&=-\sum_{i,j,k}f_{k}\Gamma_{ij}^{k}\omega^{i}\cdot\omega^{j}+\sum
_{i,j}\partial_{i}\left(  \partial_{j}\left(  f\right)  \right)
\omega
^{i}\cdot\omega^{j}\\
&  =\sum_{i,j}\left(  \partial_{i}\left(  \partial_{j}\left(
f\right) \right)  -\sum_{k}\Gamma_{ij}^{k}\partial_{k}\left(
f\right)  \right) \omega^{i}\cdot\omega^{j}\\
&  =\sum_{i,j}\left(
\partial_{i}\left(  f_{j}\right) -\sum_{k}\Gamma
_{ij}^{k}f_{k}\right) \frac{\tau_{i}\cdot\tau_{j}}{f_{i}f_{j}}.
\end{align*}
In other words,
\[
2Q_{ij}=\frac{\partial_{i}\left(  f_{j}\right) +\partial_{j}\left(
f_{i}\right)  }{f_{i}f_{j}}-\sum_{k}(\Gamma_{ij}^{k}+\Gamma_{ji}^{k}%
)\frac{f_{k}}{f_{i}f_{j}}.
\]

In what follows we shall assume that the connection $\nabla$ is
torsion-free, $\partial_{i}=\frac{\partial}{\partial x{i}}$, in
some local coordinates,  and as a result $\Gamma^i_{jk} =
\Gamma^i_{kj}$.

Summarizing, we get the following result.

\begin{theorem}
The foliation defined by level sets of a function $f$ is totally
geodesic in the torsion-free connection $\nabla$ if and only if
the function $f$ satisfies
the following system of differential equations:%
\begin{equation}
f_{j}^{2}(f_{ii}-\sum_{k}\Gamma_{ii}^{k}\
f_{k})-2f_{i}f_{j}(f_{ij}-\sum_{k}\Gamma_{ij}^{k}f_{k})+f_{i}^{2}(f_{jj}-\sum_{k}\Gamma_{jj}^{k}f_{k})=0
\label{geodcondfor f}%
\end{equation}
for all $i<j,$ where
$$
f_{ij}=(\partial_{i} (f_{j}) +\partial_{j} (f_{i}))/2.
$$
\end{theorem}

We call such a system a \emph{flex system.}

Consider a $d$-web $W_{d}$, $d \geq n + 1$, formed by the level
sets of functions $f_{\alpha}(x_{1},...,x_{n}), \alpha= 1, \dots,
d$.

The web $W_{d}$ is said to be a \emph{geodesic web} if the leaves
of all its foliations are totally geodesic.

\begin{corollary}
A web $W_{d}$ is geodesic if and only if conditions
$(\ref{geodcondfor f})$ are satisfied for all web functions
$f_{\alpha}(x_{1},...,x_{n}), \alpha= 1, \dots, d.$
\label{geodesic web}
\end{corollary}

The conditions for a planar 3-web to be geodesic were found in
another form in \cite{Ch74}.

In dimension $n=2$ we get the only one differential equation which
we have named the \textit{flex equation} (see \cite{GL06}) for the
case of flat
connections:%
\begin{equation}
f_{2}^{2}~(f_{11}-\Gamma_{11}^{1}f_{1}-\Gamma_{11}^{2}f_{2})-2f_{1}%
f_{2}~(f_{12}-\Gamma_{12}^{1}f_{1}-\Gamma_{12}^{2}f_{2})
+f_{1}^{2}~(f_{22}-\Gamma_{22}^{1}f_{1}-\Gamma_{22}^{2}f_{2})=0.
\label{geodcond n=2}%
\end{equation}

In what follows, we shall often use the following definition.

\begin{definition}
The \textit{flex} of a function $f(x,y)$ is
\[
\operatorname{Flex}\,f=f_{y}^{2}f_{xx}-2f_{x}f_{y}f_{xy}+f_{x}^{2}f_{yy}%
=-\det\left(
\begin{array}
[c]{ccc}%
f_{xx} & f_{xy} & f_{x}\\
f_{xy} & f_{yy} & f_{y}\\
f_{x} & f_{y} & 0
\end{array}
\right).
\]

\end{definition}

\section{Factorization of Flex Equation and Euler Equation}

The flex equation possesses the infinite-dimensional group of
gauge symmetries of the form
$$
f\longmapsto\Phi\left(  f\right)
$$
for any diffeomorphism $\Phi:\mathbb{R\rightarrow}\mathbb{R}$ .
Factorization of this equation with respect to this pseudogroup
leads us to the classical Euler equation as well as to some
generalizations of it.

In order to factorize the flex equation, we define the function
$$
w=\frac{f_{y}}{f_{x}},
$$
which is the first-order differential invariant of the  gauge
pseudogroup.

Indeed, we have:%
\[
\frac{\Phi\left(  f\right)  _{y}}{\Phi\left(  f\right)  _{x}}=\frac{f_{y}%
}{f_{x}},
\]
and the flex equation can be rewritten in terms of this invariant as follows:%
\[
\partial_{y}w-w\partial_{x}w=\Pi_{11}^{2}w^{3} - 3 \Pi_{12}^{2}w^{2}-3
\Pi_{12}^{1}w+\Pi_{22}^{1},
\]
where
$$
\Pi_{22}^{1} = \Gamma_{22}^{1}, \;\;
\Pi_{12}^{1}=-\frac{1}{3}(\Gamma_{22}^{2}-2\Gamma_{12}^{1}), \;\;
\Pi_{12}^{2}=-\frac{1}{3} (\Gamma_{11}^{1}-2\Gamma_{12}^{2}), \;\;
\Pi_{11}^{2}=\Gamma_{11}^{2}
$$
are the Thomas parameters (see \cite{T25}).

We call the above equation \emph{ Euler's equation associated with
the connection} $\nabla.$

Given Cauchy data $w(x,0)=w_{0}(x),$ one can solve the above Euler
equation by the standard method of characteristics, and then find
the function $f\left( x,y\right)  $ as a first integral of the
vector field:
\[
\partial_{y}  -w\partial_{x}.
\]

\section{Geodesic 4-Webs and Projective Structures}

Let us rewrite equation (\ref{geodcond n=2}) as follows:%
\begin{equation}
\Pi_{22}^{1}f_{1}^{3}-3\Pi_{12}^{1}f_{1}^{2}%
f_{2}-3\Pi_{12}^{2}f_{1}f_{2}^{2}+\Pi_{11}^{2}%
f_{2}^{3}=\operatorname{Flex}\,f. \label{geodcond2 n=2}%
\end{equation}

Note that equation (\ref{geodcond2 n=2}) appeared also in the
paper \cite{Gou87}, where the author studied the general theory of
systems of linear second-order PDEs.

We shall consider equation (\ref{geodcond2 n=2})  as a linear
equation for the components $\Pi_{jk}^{i}$ of the connection.

Remind that two affine connections, say $\nabla$ and
$\widetilde{\nabla},$ are \textit{projectively equivalent} if
there is a differential $1$-form $\rho$ such that
\[
\nabla_{X} (Y) - \widetilde{\nabla}_{X} (Y) =\rho (X) Y+\rho (Y) X
\]
for all vector fields $X$ and $Y$ (see, for example, \cite{NS94},
p. 17).

A \textit{projective structure} on a manifold may be defined as a
class of projectively equivalent affine connections.

The coefficients $\Pi^k_{ij}$ completely determine the equivalence
class of the connection given by $(\Gamma^i_{jk})$.

Assume that a geodesic $4$-web is given by web functions
\[
f_{1}(x,y),f_{2}(x,y),f_{3}(x,y),f_{4}(x,y).
\]

Then equation (\ref{geodcond2 n=2}) gives the following linear
system with respect to the Thomas parameters (or the so-called
projective connection):
\begin{equation}
\renewcommand{\arraystretch}{1.3}%
\begin{array}
[c]{l}%
\Pi_{22}^{1}f_{1,1}^{3}-3 \Pi_{12}^{1} f_{1,1}%
^{2}f_{1,2}-3 \Pi_{12}^{2} f_{1,1}f_{1,2}^{2}+\Pi_{11}^{2}f_{1,2}^{3}=\operatorname{Flex}\,f_{1},\\
\Pi_{22}^{1}f_{2,1}^{3}-3 \Pi_{12}^{1} f_{2,1}%
^{2}f_{2,2}-3 \Pi_{12}^{2} f_{2,1}f_{2,2}^{2}+\Pi_{11}^{2}f_{2,2}^{3}=\operatorname{Flex}\,f_{2},\\
\Pi_{22}^{1}f_{3,1}^{3}-3 \Pi_{12}^{1} f_{3,1}%
^{2}f_{3,2}-3 \Pi_{12}^{2} f_{3,1}f_{3,2}^{2}+\Pi_{11}^{2}f_{3,2}^{3}=\operatorname{Flex}\,f_{3},\\
\Pi_{22}^{1}f_{4,1}^{3}-3 \Pi_{12}^{1} f_{4,1}%
^{2}f_{4,2}-3 \Pi_{12}^{2}
f_{4,1}f_{4,2}^{2}+\Pi_{11}^{2}f_{4,2}^{3}=\operatorname{Flex}\,f_{4}.
\end{array}
\renewcommand{\arraystretch}{1} \label{system for gammas}%
\end{equation}

Solving system (\ref{system for gammas}), we get%
\begin{equation}
\renewcommand{\arraystretch}{4}%
\begin{array}
[c]{ll}%
\Pi_{22}^{1}=\displaystyle\sum_{i=1}^{4}\displaystyle\frac{\displaystyle\prod_{k\neq i}f_{k,2}}%
{\displaystyle\prod_{k\neq i}\,J(f_{i},f_{k})}\operatorname{Flex}\,f_{i}, & \\
\Pi_{11}^{2}=\displaystyle\sum_{i=1}^{4}\displaystyle\frac{\displaystyle\prod_{k\neq i}f_{k,1}}%
{\displaystyle\prod_{k\neq
i}\,J(f_{i},f_{k})}\operatorname{Flex}\,f_{i}; &
\end{array}
\renewcommand{\arraystretch}{1} \label{solution1 of linear system}%
\end{equation}

\begin{equation}
\renewcommand{\arraystretch}{4}%
\begin{array}
[c]{ll}%
-3\Pi_{12}^{1}=\displaystyle\sum_{i=1}^{4}\displaystyle\frac{\displaystyle\sum_{k\neq
i}f_{k,1}\displaystyle\prod_{l\neq i,k}f_{l,2}}{\displaystyle\prod_{m\neq i}\,J(f_{i},f_{m}%
)}\operatorname{Flex}\,f_{i}, & \\
-3\Pi_{12}^{2}=\displaystyle\sum_{i=1}^{4}\displaystyle\frac{\displaystyle\sum_{k\neq
i}f_{k,2}\displaystyle\prod_{l\neq i,k}f_{l,1}}{\displaystyle\prod_{m\neq i}\,J(f_{i},f_{m}%
)}\operatorname{Flex}\,f_{i}, &
\end{array}
\renewcommand{\arraystretch}{1} \label{solution2 of linear system}%
\end{equation}
where
\[
\,J(f_{i},f_{j})=\det\left(
\begin{array}
[c]{cc}%
f_{i,1} & f_{i,2}\\
f_{j,1} & f_{j,2}%
\end{array}
\right)
\]
is the Jacobian of the functions $f_{i}(x,y)$ and $f_{j}(x,y).$
Note that in the planar case we have $\Pi_{22}^{2} = -
\Pi_{12}^{1}$ and $\Pi_{11}^{1} = - \Pi_{12}^{2}$. Thus, formulas
(\ref{solution1 of linear system}) and (\ref{solution2 of linear
system}) give all Thomas parameters $\Pi^i_{jk}$.

Remark that system (\ref{system for gammas}) is invariant with
respect to the gauge transformations $f_i \rightarrow \Phi_i
(f_i), \, i = 1,2,3,4$, and therefore, solutions (\ref{solution1
of linear system}) and (\ref{solution2 of linear system}) do not
depend on functions $\{f_i\}$, but they are completely determined
by the geodesic 4-web.

Summarizing, we get the following result.

\begin{theorem}
Any planar $4$-web defines a unique projective structure in the
plane in such a way that the leaves of the foliations are
geodesics of this projective structure.
\end{theorem}

This theorem allows us to get the following Gronwall-type theorem
(see \cite{Gr12}).

\begin{corollary} Any diffeomorphism sending a planar $4$-web
$W_4$ into a planar $4$-web $\widetilde{W}_4$ is a projective
transformation of the corresponding projective structures.
\end{corollary}

\section{Geodesic Webs and Symmetric Projective Structures}

We say that a projective structure is \textit{symmetric} if the
class of projectively equivalent affine connections contains an
affine symmetric connection.

Let us consider the case when the projective structure determined
by a 4-web is symmetric. Then $\nabla R = 0$ for some affine
connection from the class.

Denote
\begin{equation}
\Gamma_{12}^{1}=\sigma,
\;\alpha=\Gamma_{22}^{2}-2\Gamma_{12}^{1},\;
\beta=\Gamma_{11}^{1}-2\Gamma_{12}^{2}, \;\Gamma_{12}^{2}=\tau.
\label{notation for gammas}%
\end{equation}
Then
\begin{equation}
\Gamma_{11}^{1}=2\tau+\beta, \;\Gamma_{22}^{2}=2\sigma+\alpha.
\label{gammai-ii}%
\end{equation}

In order to simplify formulae, we choose coordinates $x, y$ in the
plane in such a way that $f_{1} (x,y)  =x,\, f_{2}(x,y) =y.$

Then (\ref{solution1 of linear system}) gives
\begin{equation}
\Gamma_{11}^{2}=0, \;\Gamma_{22}^{1}=0, \label{gamma=0}
\end{equation}
and (\ref{solution2 of linear system}) becomes
\begin{equation}
\alpha =\frac{f_{4,y}\operatorname{Flex}\,f_3}{f_{3,x}f_{3,y}
\Delta}- \frac{f_{3,y}\operatorname{Flex}\,f_4}{f_{4,x}f_{4,y}
\Delta},\;\; \beta
=-\frac{f_{4,x}\operatorname{Flex}\,f_3}{f_{3,x}f_{3,y} \Delta}
+\frac{f_{3,x}\operatorname{Flex}\,f_4}{f_{4,x}f_{4,y} \Delta},
\label{alpha-beta}
\end{equation}
where
\[
\Delta=f_{3,x}f_{4,y}-f_{3,y}f_{4,x}.
\]

By a straightforward computation, one can find the components of
the curvature tensor $R\left(  \nabla\right) =\{R_{jkl}^{i}\}$ of
the connection :%
\[
R\left(  \nabla\right)  =%
\begin{Vmatrix}
R_{112}^{1} & R_{212}^{1}\\
R_{112}^{2} & R_{212}^{2}%
\end{Vmatrix}
=%
\begin{Vmatrix}
\sigma_{x}-2\tau_{y}-\beta_{y}+\sigma\tau & -\sigma_{y}+\sigma^{2}%
+\sigma\alpha\\
\tau_{x}-\tau^{2}-\tau\beta &
2\sigma_{x}-\tau_{y}+\alpha_{x}-\sigma\tau
\end{Vmatrix}
.
\]

Computing components of $\nabla R$, we get the following system of
differential equations for components $\tau$ and $\sigma$:
\begin{equation}
\renewcommand{\arraystretch}{1.3}
\begin{array}{ll}
{\sigma_{xx}}-2\,{\tau_{xy}}-{\beta_{xy}}  &  =(2\tau\,+\beta){\sigma_{x}%
}-2\,\sigma{\tau}_{x}-(3\,\tau+\beta)(2{\tau_{y}}+{\beta_{y})}+2\,\sigma
\tau (2\tau +\beta),\\
\sigma_{xy}-2\,\tau_{yy}-\beta_{yy}  &  =(3\,\sigma+\alpha)\sigma
_{x}-(7\,\sigma+2\alpha)\tau_{y}-(3\,\sigma+\alpha)\beta_{y}-2\,\tau\sigma
_{y}+2\,\sigma\tau (2\, \sigma+\alpha),\\
\alpha_{xx}+2\,\sigma_{xx}-\tau_{xy}  &  =(3\,\tau+\beta)\alpha_{x}%
+(7\tau+2\beta)\,\sigma_{x}+2\,\sigma\tau_{x}-(3\,\tau+\beta)\tau
_{y}-2\,\sigma\tau (2\, \tau+\beta),\\
\alpha_{xy}+2\,\sigma_{xy}-\tau_{yy}  &  =(3\,\sigma+\alpha)\alpha
_{x}+2(3\,\sigma+\alpha)\sigma_{x}-(2\,\sigma+\alpha)\tau_{y}+2\,\tau
\sigma_{y}-2\,\sigma\tau (2\, \sigma+\alpha),
\end{array} \label{sys1}
\renewcommand{\arraystretch}{1}
\end{equation}
and%
\begin{equation}
\renewcommand{\arraystretch}{1.3}
\begin{array}{ll}
\tau_{xx}  & =2(3\,\tau+\beta)\tau_{x}+\tau\beta_{x}-2\,\tau
(2\,\tau^{2}+3\,\tau \beta+ \beta^{2}),\\
\tau_{xy}  &  =\tau\sigma_{x}+\tau\alpha_{x}+(3\,\tau+\beta)\tau_{y}%
+2\,\tau\beta_{y}+2\,\sigma\tau_{x}-2\,\sigma\tau (2\, \tau+\beta),\\
\sigma_{xy}  &
=(3\,\sigma+\alpha)\sigma_{x}+2\,\sigma\alpha_{x}+\sigma
\tau_{y}+2\,\tau\sigma_{y}+\sigma\beta_{y}-2\,\sigma\tau (2\, \sigma+\alpha),\\
\sigma_{yy}  &  =3(2\sigma+\alpha)\sigma_{y}-2\,\sigma\alpha^{2}%
+3\alpha\,\sigma_{y}+\sigma\alpha_{y}-2\,\sigma^{2}
(2\,\sigma+3\,\alpha).
\end{array} \label{sys2}
\renewcommand{\arraystretch}{1}
\end{equation}

Consider a system consisting of the first and the last equations
of (\ref{sys1}) and (\ref{sys2}). Solving this system, we find all
second derivatives $\tau_{xx},\tau_{xy},\tau_{yy}$ and
$\sigma_{xx},\sigma _{xy},\sigma_{yy}$ in terms of $\alpha, \beta,
\sigma, \tau$ and first-order derivatives of $\sigma$
and $\tau:$%
\begin{equation}
\renewcommand{\arraystretch}{1.3}
\begin{array}{ll}
\sigma_{xx}  & =2\,\sigma\tau_{x}+(4\,\tau
+\beta)\sigma_{x}+(\tau-\beta)\beta_{y}+2\,\tau\alpha_{x}+\beta_{xy}-2\,\sigma\tau(
2\,\tau+\beta),\\
\sigma_{xy}  & =(3\,\sigma + \alpha)\sigma_{x}+2\,\sigma\alpha_{x}+\sigma\tau_{y}%
+2\,\tau\sigma_{y}+\sigma\beta_{y}-2\,\sigma\tau(2\sigma +\alpha),\\
\sigma_{yy}  & =3(2\sigma + \alpha)\sigma_{y} +
\sigma\alpha_{y}-2\,\sigma (\alpha^{2}+2\,\sigma^{2}+3\,\sigma \alpha),\\
\tau_{xx}  & =3(2\,\tau + \beta)\tau_{x} +\tau\beta_{x} -2\,\tau(\beta^{2}+2\,\tau^{2}+3\,\tau \beta),\\
\tau_{xy}  & =\tau \sigma_{x}+\tau\alpha_{x}+(3\,\tau + \beta)
\tau_{y}+2(\tau\beta _{y}+\,\sigma\tau_{x})-2\,\sigma\tau
(2\,\tau + \beta),\\
\tau_{yy}  & =(\sigma - \alpha)\alpha_{x} + (4\,\sigma +
\alpha)\tau_{y} +
2(\tau\sigma_{y}+\,\sigma\beta_{y})+\alpha_{xy}-2\,\sigma \tau
(2\, \sigma + \alpha).
\end{array} \label{sys3}
\renewcommand{\arraystretch}{1}
\end{equation}
Substituting these expressions into the second and third equations
of (\ref{sys1}), we get the following relations for $\alpha$ and
$\beta$:
\begin{equation}
\renewcommand{\arraystretch}{1.3}
\begin{array}{ll}
\alpha_{xx}+2\beta_{xy} &  =\beta\alpha_{x}+2\beta\beta_{y},\\
2\alpha_{xy}+\beta_{yy} & =2\alpha\alpha_{x}+\alpha\beta_{y}.
\end{array} \label{sys4}
\renewcommand{\arraystretch}{1}
\end{equation}
Remark that for affine symmetric connections we have tr$\;R=0,$
and it is easy to check that this is the only compatibility
condition for system (\ref{sys3}) under conditions (\ref{sys4})
for $\alpha$ and $\beta.$

Moreover, adding the condition $\operatorname*{tr}\;R=0,$ or
\begin{equation}
\alpha_{x}-\beta_{y}+3(\sigma_{x}-\tau_{y})=0,\label{sys5}%
\end{equation}
does not produce any new compatibility condition.

In other words, assuming conditions (\ref{sys4}) for $\alpha$ and
$\beta,$ the PDE system of equations  (\ref{sys3}) and
(\ref{sys5}) with respect to  $\sigma$ and  $\tau$ is a formally
integrable system of finite type. It is easy to see that the
solution space of the system has dimension 5.

Summarizing, we arrive at the following result.

\begin{theorem}
\begin{description}
\item[\rm{(i)}]  The projective structure defined by a planar
$4$-web is symmetric if and only if conditions $(\ref{sys4})$
hold.

\item[\rm{(ii)}]  A planar $4$-web given by functions $\{x, y, f_3
(x, y), f_4 (x, y)\}$ is locally equivalent to a geodesic $4$-web
on an affine symmetric surface if and only if conditions
$(\ref{sys4})$ hold for the functions $\alpha$ and $\beta$ defined
by equation $(\ref{alpha-beta})$.
\end{description}
\end{theorem}

\begin{remark}
If a planar $4$-web satisfies conditions of the above theorem,
then there exists a $5$-dimensional family
$\nabla_{t},~t\in\mathbb{R}^{5},$ of affine symmetric connections
such that all leaves of the web are geodesics for $\nabla_{t}.$
Take one of them, say $\nabla_{t_{0}}.$ Then for cases $R
(\nabla_{t_{0}}) \neq 0,$ if $\det R (\nabla_{t_{0}}) > 0$, then
the $4$-web is equivalent to a geodesic $4$-web either on the
standard $2$-sphere or to a geodesic $4$-web on the Lobachevskii
plane. If $\det R (\nabla_{t_{0}}) < 0$, then the $4$-web is
equivalent to a geodesic $4$-web on the de Sitter plane. If $\det
R (\nabla_{t_{0}}) = 0$, then the $4$-web is equivalent to a
geodesic $4$-web on the affine torus or on the Klein bottle
$($see, for example, \rm{\cite{F76}}$)$.
\end{remark}

\section{Planar Geodesic \boldmath${d}$-Webs and Projective Structures}

Consider now a planar $d$-web $W_{d}$ defined by $d$ web functions
$f_{\alpha}(x,y),\;\alpha=1,...,d.$ Such a web has $\binom{d}{4}$
$4$-subwebs $[\alpha,\beta,\gamma,\varepsilon]$ defined by the
foliations $X_{\alpha},X_{\beta},$ \newline $X_{\gamma},$ and
$X_{\varepsilon}, \; \alpha,\beta,\gamma,\varepsilon=1,...,d.$

If a $d$-web $W_d$ is geodesic, then each of $\binom{d}{4}$ its
$4$-subwebs $[\alpha,\beta,\gamma,\varepsilon]$ is also geodesic,
and by Theorem 8, each of them determines its own unique
projective structure. These $d-4$ projective structures coincide
with a projective structure defined by one of them, let's say, by
the 4-subweb $[1,2,3,4]$, if $d-4$ second-order invariants given
by the flex equations vanish.

Thus we have proved the following result.

\begin{theorem}
\begin{description}
\item[\rm{(i)}] A planar $d$-web, $d \geq 5$, defined by  web
functions $<f_1, \dots, f_d>$ is geodesic if and only if the
functions $f_5, \dots, f_d$ satisfy flex equations
$(\ref{geodcond2 n=2})$, in which the components $\Pi^1_{22},
\Pi^2_{11}, \Pi^1_{12}$ and $\Pi^2_{12}$ are given by formulae
$(\ref{solution1 of linear system})$ and $(\ref{solution2 of
linear system})$.

\item[\rm{(ii)}] A planar $d$-web, $d \geq 5$, defined by web
functions $<x, y, f_3, \dots, f_d>$ is locally equivalent to a
geodesic $d$-web with respect to a symmetric projective structure
if and only if the $d-4$ conditions for geodesicity mentioned
above are satisfied and, in addition, conditions $(\ref{sys4})$
hold.
\end{description}
\end{theorem}

\begin{remark}
Theorem 4 shows that one has $d - 4$ second-order conditions on
web functions and two more fourth-order conditions to have a
geodesic $d$-web in a symmetric projective structure.
\end{remark}

\section{Planar Linear Webs}

In what follows, we shall use  coordinates $x, y$ on the plane in
which the Christoffel symbols $\Gamma_{jk}^{i}$  vanish. We shall
assume that a $d$-web
$W_{d}$\ is formed by the level sets of  web functions $f_{1}%
(x,y),\;f_{2}(x,y),.....,$ $f_{d}(x,y).$

The following theorem, which immediately follows from formula
(\ref{geodcond2 n=2}), gives a criterion for $W_{d}$ to be  linear
in the coordinates $x,y$.

\begin{theorem}
The $d$-web $W_{d}$ is a linear if and only if the web functions
are solutions of the differential equation
\begin{equation}
\operatorname{Flex}\,f=0.\label{lincondition}%
\end{equation}

\end{theorem}

Note that in algebraic geometry the linearity condition $\operatorname{Flex}%
\,f=0$ is also the necessary and sufficient condition for a point
$(x,y)$ to be a flex of the curve defined by the equation
$f(x,y)=0$. Here, in (\ref{lincondition}),
$\operatorname{Flex}\,f=0$ is the equation for finding the
function $f(x,y)$ (it should be satisfied for all points $(x,y)$).
This is a reason that we call equation (\ref{lincondition}) the
\emph{flex equation}.

We shall show  how to integrate flex equation
(\ref{lincondition}). The main idea of integration is that the
factorization of the flex-equation $\operatorname{Flex}\,f=0$ with
respect to the diffeomorphism group produces an Euler equation.

Namely, as we have seen, passing to the differential invariant
$w=\frac{f_{x}}{f_{y}}$ allows us to reduce the order of the
equation. Let us rewrite the flex equation in the form
\begin{equation}
\partial_{x}\left(  \frac{f_{x}}{f_{y}}\right)  -\left(  \frac{f_{x}}{f_{y}%
}\right)  \partial_{y}\left(  \frac{f_{x}}{f_{y}}\right)  =0.\label{Euler1}%
\end{equation}
Then integration of the flex equation is equivalent to solution of
the
following system:%
\[
\left\{
\begin{array}
[c]{ll}%
\partial_{x}w-w\partial_{y}w & =0,\\
\partial_{x}f-w\partial_{y}f & =0.
\end{array}
\right.
\]

The first equation of the system
\[
\partial_{x}w-w\partial_{y}w=0
\]
is the classical Euler equation in gas-dynamics.

Solutions of this equation are well-known. Namely, if $w_{0}\left(
y\right) =\left.  w\right\vert _{x=0}$ is the Cauchy data, then
the solution $w(x,y)$ can be found from the system
\begin{equation}
\left\{
\begin{array}
[c]{lll}%
y+w_{0}\left(  \lambda\right)  x-\lambda & = & 0,\\
w(x,y)-w_{0}\left(  \lambda\right)   & = & 0
\end{array}
\right.  \label{EulerWeb}%
\end{equation}
by elimination of the parameter $\lambda.$

Further, if $w$ is a solution of the Euler equation, then the
functions $w$
and $f$ are both first integrals of the vector field $\partial_{x}%
-w\partial_{y},$ and therefore, $f=\Phi(w)$ for some smooth
function $\Phi$.

Summarizing, we get the following description of web functions of
linear webs.

\begin{proposition}
\label{3weband Euler1}The web functions
$f_{1}(x,y),f_{2}(x,y),...,$ $f_{d}(x,y)$ of a linear $d$-web have
the form
\[
f_{1}(x,y)=\Phi_{1}(w_{1}(x,y)),\;f_{2}(x,y)=\Phi_{2}(w_{2}(x,y)),...,\;f_{d}%
(x,y)=\Phi_{d}(w_{d}(x,y)),
\]
where $w_{1}(x,y),w_{2}(x,y),...,$ $w_{d}(x,y)$ are  distinct
solutions of the Euler equation, and
$\Phi_{1},\Phi_{2},...,\Phi_{d}$ are smooth functions.
\end{proposition}

In particular, using the gauge transformations, we can take
\[
f_{1}(x,y)=w_{1}(x,y),\;f_{2}(x,y)=w_{2}(x,y),....,\;f_{d}(x,y)=w_{d}(x,y).
\]
Therefore, the above proposition yields the following description
of web functions for linear $d$-webs.

\begin{theorem}
\label{3weband Euler2} Web functions of linear $d$-webs can be
chosen as $d$ distinct solutions of the Euler equation.
\end{theorem}

\begin{example}
\label{exparabola} Assume that for a linear 3-web  we have
$f_{1}(x,y)=x,$ $f_{2}(x,y)=y\;$and$\;f_{3}(x,y)=f(x,y).$ Taking
$w_{0}\left( y\right)  =-2\sqrt{-y},$ we get from
$(\ref{EulerWeb})$ that $w=-2(x+\sqrt {x^{2}-y})$ and
$f=x+\sqrt{x^{2}-y}$. The leaves of the third foliation are the
tangents to the parabola $y=x^{2}.$
\end{example}

\begin{example}
Assume that for a linear 5-web we have $f_{1}(x,y)=x$ and
$f_{2}(x,y)=y.$ Taking $(w_{3})_{0} (y) = -2\sqrt{-y},
\;(w_{4})_{0} (y) = y$ and $(w_{5})_{0} (y) = 2y,$ we get the
linear $5$-web with remaining three web functions
\[
f_{3}=x+\sqrt{x^{2}-y}, \; f_{4}=\displaystyle \frac{y+1}{1-x},
\;\; f_{5}= \displaystyle                   \frac{y}{1-2x}.
\]
The last three foliations of this $5$-web are the tangents to the
parabola $y=x^{2}$ $($see Example $\ref{exparabola})$ and the
straight lines of the pencils with the centers $(1,-1)$ and
$(\frac{1}{2},0).$
\end{example}

\section{ Geodesic Webs on Surfaces of Constant Curvature}

\begin{theorem}
Let $(M,g)$ be a surface of constant curvature with the metric
tensor
\begin{equation}
g=\frac{dx^{2}+dy^{2}}{(1+\kappa(x^{2}+y^{2}))^{2}}%
,\label{metric of constant curv}%
\end{equation}
where $\kappa$ is a constant. Then the level sets of a function
$f(x,y)$ are geodesics of the metric if and only if the function
$f$ satisfies the flex equation
\begin{equation}
\operatorname{Flex}\;f=\frac{2\kappa\left(  xf_{x}+yf_{y}\right)
\left( f_{x}^{2}+f_{y}^{2}\right)  }{1+\kappa\left(
x^{2}+y^{2}\right)
}.\label{flex equation}%
\end{equation}

\end{theorem}

\begin{proof}
It is easy to see that the Christoffel symbols $\Gamma_{jk}^{i}$
are the following
\[
\renewcommand{\arraystretch}{1.3}%
\begin{array}
[c]{l}%
\Gamma_{12}^{1}=\Gamma_{22}^{2}=\Gamma_{21}^{1}=-\Gamma_{11}^{2}=-2\kappa
x_{2}b,\;\\
\Gamma_{11}^{1}=\Gamma_{21}^{2}=\Gamma_{12}^{2}=-\Gamma_{22}^{1}=-2\kappa
x_{1}b,\;\;
\end{array}
\renewcommand{\arraystretch}{1}
\]
where
\[
b=\frac{1}{1+\kappa\left(  x^{2}+y^{2}\right)  }.
\]

Substituting these values into the right-hand side of formula
(\ref{geodcondfor f}), we get formula (\ref{flex equation}).
\end{proof}

Remark that functions of the form $f(x,y)=\Phi\left(
\frac{y}{x}\right)  $ are solutions of the flex equation, and
therefore  the level sets of these functions are geodesic.

As we have seen (see Theorem 4), a geodesic $d$-web uniquely
defines a projective structure  provided that  $d-4$ additional
second-order invariants vanish. This leads us to the following
Gronwall--type theorem (cf. Corollary 4):

\begin{theorem}
Suppose that $W_{d}$, $d \geq 4,$ is a geodesic $d$-web given on a
surface $(M,g)$ of constant curvature for which  $d-4$ additional
second-order invariants vanish. Then any mapping of $W_{d}$ on a
geodesic web $\widetilde{W}_{d}$ is a projective transformation.
\end{theorem}

\begin{proof}
We give an alternative proof. By Beltrami's theorem (see
\cite{Be65}), the surface $(M,g)$ can be mapped by a
transformation $\phi$ onto a plane, and the mapping $\phi$ sends
the geodesics of $(M, g)$ into straight lines. Thus the mapping
$\phi$ linearizes the webs $W_{d}$ and $\widetilde{W}_{d}.$ In
\cite{GL06} (see Remark on p. 99) it was proved that for $d \geq
4$ the Gronwall conjecture is valid, i.e., there exists a unique
transformation sending $W_{d}$ and $\widetilde{W}_{d}$ into linear
$d$-webs, and this transformation is projective. The mapping of
$W_{d}$ onto $\widetilde{W}_{d}$ induces a transformation of the
corresponding linear $d$-webs. The latter transformation is
projective. As a result, the mapping of $W_{d}$ onto
$\widetilde{W}_{d}$ is also projective.
\end{proof}

\section{Geodesic Webs on Surfaces in $\mathbb{R}^3$}

\begin{proposition}
Let $(M,g)\subset\mathbb{R}^{3}$ be a surface defined by an
equation $z=z(x,y)$ with the induced metric $g$ and the
Levi-Civita connection
$\nabla.$ Then the flex equation takes the form%
\begin{equation}
\operatorname{Flex}\;f=\frac{z_{x}f_{x}+z_{y}f_{y}}{1+z_{x}^{2}+z_{y}^{2}%
}(f_{y}^{2}z_{xx}-2f_{x}f_{y}z_{xy}+f_{x}^{2}z_{yy}). \label{flex of z(x,y)}%
\end{equation}

\end{proposition}

\begin{proof}
To prove formula (\ref{flex of z(x,y)}), note that the metric
induced on a surface $z=z(x,y)$ is
\[
g=ds^{2}=(1+z_{x}^{2})dx^{2}+(1+z_{y}^{2})dy^{2}+2z_{x}z_{y}dxdy.
\]
Computing the Christoffel symbols, we get that%
\[
\renewcommand{\arraystretch}{2}%
\begin{array}
[c]{c}%
\Gamma_{11}^{1}=\displaystyle\frac{z_{x}z_{xx}}{1+z_{x}^{2}+z_{y}^{2}%
},\;\Gamma_{12}^{1}=\Gamma_{21}^{1}=\displaystyle\frac{z_{x}z_{xy}}%
{1+z_{x}^{2}+z_{y}^{2}},\;\Gamma_{22}^{1}=\displaystyle\frac{z_{x}z_{yy}%
}{1+z_{x}^{2}+z_{y}^{2}},\\
\Gamma_{11}^{2}=\displaystyle\frac{z_{y}z_{xx}}{1+z_{x}^{2}+z_{y}^{2}%
},\;\Gamma_{12}^{2}=\Gamma_{21}^{2}=\displaystyle\frac{z_{y}z_{xy}}%
{1+z_{x}^{2}+z_{y}^{2}},\;\Gamma_{22}^{2}=\displaystyle\frac{z_{x}z_{yy}%
}{1+z_{x}^{2}+z_{y}^{2}}.
\end{array}
\renewcommand{\arraystretch}{1}
\]
Applying these formulas to the right-hand side of (\ref{geodcond
n=2}), we get formula (\ref{flex of z(x,y)}).
\end{proof}

For example, functions of the form $f(x,y)=\Phi\left(
\frac{x}{y}\right)  $ are solutions of the flex equation if and
only if the function $z\left(
x,y\right)  $ satisfies one of the following equations:%
\[
yz_{x}-xz_{y}=0,
\]
or%
\[
x^2 z_{xx}+2xy z_{xy}+y^2 z_{yy}=0.
\]
The general  solution of the first equation have the form $z =
\Omega (x^{2} + y^{2}),$ and of the second one $z(x, y) =x \Psi
(\frac{y}{x}) +\Theta (\frac{y}{x}) $, where $\Omega (\alpha),
\Psi (\alpha)$ and $\Theta(\alpha)$ are arbitrary smooth
functions.

If we assume that the foliations $\{x=\operatorname*{const}.\}$
and $\{ y=\operatorname*{const}.\}$ are geodesic on the surface
$z=z (x,y)$, then the flex equation gives $z_{xx}=z_{yy}=0,$ and
therefore $z=axy+bx+cy+d$ for some constants $a, b, c$, and $d$,
and formula (\ref{flex for z(x,y)}) takes the form
$$
\operatorname{Flex}\;f=-\frac{2a f_{x} f_{y}((ay+b)f_{x}+(ax+c)f_{y})}{1+(ay+b)^{2}+(ax+c)^{2}%
}. \label{flex for z(x,y)}%
$$

{\emph{Authors' addresses:} }

{Deparment of Mathematical Sciences, New Jersey Institute of Technology,
University Heights, Newark, NJ 07102, USA; vladislav.goldberg@gmail.com}

{Department of Mathematics, The University of Tromso, N9037, Tromso, Norway;
lychagin@math.uit.no}

\end{document}